\documentclass[11pt,a4paper]{article}

\usepackage{amsmath}
\usepackage{amsthm}
\usepackage{amsfonts}
\usepackage[koi8-r]{inputenc}

\usepackage[russian]{babel}

\parskip\smallskipamount
\theoremstyle{remark}

\renewcommand{\phi}{\varphi}
\renewcommand{\epsilon}{\varepsilon}

\parskip \smallskipamount
\parskip \smallskipamount

\title{{\Large \bf Some open problems related to stability}}

\author{S. Foss\\
Heriot-Watt University, Edinburgh and Sobolev's Institute of 
Mathematics, Novosibirsk}
\date{}

\begin{document}

\maketitle

\begin{quotation}\small
I will speak about a number of open problems in queueing. Some of them
are known for decades, some are more recent. They relate to stability and
to rare events.

There is an idea to prepare a special issue of QUESTA on open problems, and
this text may be considered as a prospective contribution to that.  The choice of open problems
reflects the speaker's own interests, and should not be taken as suggesting that these are the only,
 or even most important, problems!
\end{quotation}


\section{Multi-server queue with First-Come-First-Served discipline}

A system with a finite number of identical servers and with FCFS 
service discipline is one of the simplest models in queueing theory.
It has been known for a long time. To the best of my knowledge, Kiefer
and Wolfowitz were the first who studied it rigorously.

\subsection{Convergence in the total variation to the stationarity} 

Consider first a single-server first-come-first-served (FCFS) queue $G/G/1$ with 
interarrival times $\{t_n\}$ between the arrivals of the $n^{th}$ and $(n+1)^{st}$ customers, and service
times of the $n^{th}$ customer $\{\sigma_n\}$. Assume that the two-dimensional sequence $\{t_n,\sigma_n\}$ is
stationary ergodic and that the queue is stable, that is 
${\mathbf E} \sigma_1 < {\mathbf E} t_1$. 

Let $W_n$, $n\ge 1$ be the waiting time of customer $n$ in the system (before the start
of its service). Then
$$
W_1=x\ge 0, \quad \mbox{and} \quad W_{n+1} = \max (0, W_n+\sigma_n-t_n) \equiv
(W_n+\xi_n)^+, \quad n\ge 1
$$
where $x$ is the initial delay, $\xi_n=\sigma_n-t_n$, and $x^+ = \max (0,x)$.

Let $S_0=0$ and $S_n =\sum_1^n \xi_i$, $n\ge 1$. 
Clearly, 
$$
W_n = \max (0,x+S_n, S_n-S_1,S_n-S_2,\ldots, S_n-S_{n-1})
$$ 
which implies that there exists a proper limiting distribution
which does not depend on $x$. In other words, there
exists a unique stationary distribution of the waiting time and, for any initial delay,
there is a convergence to stationarity, and the convergence is in the total variation
norm.

There are many known ways to establish this result. The best one seems to use the
``Loynes scheme''. Without loss of generality, we may assume $(\sigma_n,t_n)$
to be defined for all $-\infty < n < \infty$. Let $\widetilde{S}_0=0$ and
$\widetilde{S}_n = \sum_{j=1}^n \xi_{-j}$,
$n\ge 1$. Then 
\begin{equation}\label{s1}
W_n =_{st} M_{n}^{(x)} := \max (0,x+\widetilde{S}_n, \widetilde{S}_{n-1}, \ldots, \widetilde{S}_1).
\end{equation}
Denote $M = \sup_{n\ge 0} \widetilde{S}_n$. 
{From} the SLLN, $\widetilde{S}_n\to -\infty$ a.s., so $M<\infty$ a.s. 
Moreover, 
the time
$$
\nu \equiv \nu^{(x)} = \max \{n\ge 0 \ : \ x+\widetilde{S}_n\ge 0\}
$$
is finite a.s.  and, therefore,
$$
M_{n}^{(x)} = M \quad \mbox{for all} \quad n> \nu^{(x)}.
$$
So, $M$ is the unique limiting distribution and
 the convergence in total variation follows from the coupling inequality:
for any $x\ge 0$, 
$$
\sup_A |{\mathbf P}(W_n\in A) - {\mathbf P} (M\in A)|
= \sup_A |{\mathbf P} (M_{n}^{(x)}\in A) - {\mathbf P}(M\in A)| 
\leq {\mathbf P} (\nu^{(x)}>n)\to 0,\quad n\to\infty.
$$ 
In particular, if $x=0$, then -- as follows from (\ref{s1}) --
the sequence $M_n :=M_n^{(0)}$, $n\ge 1$ is monotone
increasing a.s. and couples with $M$, starting from time $\nu^{(0)}+1$.

Now let $m>1$ be a positive integer and consider the $G/G/m$ FCFS queue with 
interarrival times $\{t_n\}$ and service
times $\{\sigma_n\}$. Assume again that the two-dimensional sequence $\{t_n,\sigma_n\}$ is
stationary ergodic and that the system is stable. Here the stability means
${\mathbf E} \sigma_1 < m {\mathbf E} t_1$.

Consider Kiefer-Wolfowitz vectors of virtual waiting times
${\mathbf W}_n = (W_{n1},\ldots,W_{nm})$ which satisfy the recursion 
$$
{\mathbf W}_1 = {\mathbf x}\ge {\mathbf 0} \quad \mbox{and}
\quad
{\mathbf W}_{n+1} = R ({\mathbf W}_n + {\mathbf e}_1 \sigma_n -
{\mathbf 1} t_n)^+,\quad n\ge 1
$$
where ${\mathbf x}$ is a vector of initial delays,
${\mathbf e}_1=(1,0,\ldots,0)$ is a unit vector, ${\mathbf 1}=(1,1,\ldots,1)$
is a vector of units, and operator $R$ rearranges coordinates of a vector in
 weak ascending order. In particular, $W_{n,1}$ is a waiting time of customer $n$.

It is known that, in general, in a multi-server queue there may be many stationary regimes
[\ref{Loy}]; there exist the minimal and the
maximal stationary distribution [{\ref{Bra}, \ref{BFL}], and there are 
some relations between stationary distributions [\ref{F83}]. 
Similarly to the one-dimensional case, one can define vectors 
${\mathbf M}_n^{({\mathbf x})}$
which satisfy a recursion which is more complex than (\ref{s1}). In particular, 
if there are no
initial delays, ${\mathbf x}={\mathbf 0}$, then the vectors ${\mathbf M}_n = 
{\mathbf M}_n^{({\mathbf 0})}$ are monotone weakly ascending (coordinate-wise and a.s.)
and converge a.s. to a limiting random vector which has the minimal stationary
distribution $\pi_{min}$ (see e.g. [\ref{Sto}]). But this implies only the weak convergence
of the distributions of these vectors, and not the convergence in the total variation. 

If the sequence $\{ (\sigma_n,t_n)\}$ satisfies in addition some good 
``mixing'' properties (say, is i.i.d. or regenerative), then one can again show the 
uniqueness of the stationary regime and the convergence in the total variation starting
from each initial value, using
either Harris properties (in the Markovian case, see e.g. [\ref{MT}]) 
or the renovation techniques (in a more general setting, see e.g. [\ref{Bor}]).

The {\bf conjecture} is: assuming only that the sequence $\{(\sigma_n,t_n)\}$ is stationary
ergodic, that the stability condition holds, and that the initial value is $\mathbf 0$,
then there is no need for any further restriction to establish the convergence in the
total variation:
$$
\sup_A |{\mathbf P} ({\mathbf W}_n \in A) - \pi (A)|\to 0, \quad n\to\infty.
$$

There have been a number of unsuccessful attempts to prove this conjecture
(see, e.g., [\ref{Naka}]).

\subsection{Existence of Moments}

Assume now that $\{\sigma_n\}$ and $\{t_n\}$ are two i.i.d. sequences
that do not depend on each other. Continue to assume the stability condition
$\rho:= {\mathbf E} \sigma_1/{\mathbf E} t_1 < m$ holds. Recall that in this case the
stationary distribution is unique.

Denote by $D$ the stationary waiting time in the multi-server queue.
Fix $\gamma>0$ and formulate the following question:
what are the conditions for ${\mathbf E} D^{\gamma}$ to be finite.

A correct (but partial!) answer to this question
has been obtained recently by Scheller-Wolf and Vesilo [\ref{SW-V}].  
To be completely exact, the result below was not formulated by
these authors but may be deduced from their results.

Denote by $B_I$ the integrated service time distribution,
$$
B_I(x) = 1 - \min \left( 1,\int_x^{\infty} {\mathbf P} 
(\sigma_1>y) dy \right).
$$
Let $\sigma_{I,1}, \ldots ,\sigma_{I,m}$ be i.i.d. random variables
with common distribution $B_I$.

{\bf Proposition}
Assume that $\rho$ is not an integer and denote by $k\in \{0,1,\ldots,
m-1\}$ its integer part. Then 
$$
{\mathbf E} D^{\gamma}<\infty 
\quad
\mbox{iff}
\quad
{\mathbf E} \left(\min (\sigma_{I,1},\ldots,\sigma_{I,m-k})
\right)^{\gamma} < \infty.
$$

The proofs of the results in [{\ref{SW-V}] are based
on the construction of an auxiliary, so-called ``semi-cyclic'' service
discipline. A direct proof of the proposition may be found 
in [\ref{FKor}] and is based on ideas
close to ideas of Keifer and Wolfowitz [\ref{KW}].

An {\bf open problem} here is: what are the conditions for existence
(finiteness)
of power moments of $D$ if $\rho$ is an integer.

\subsection{Rare events}

Assume again that $\{\sigma_n\}$ and $\{t_n\}$ are two i.i.d. sequences
that do not depend on each other. Continue to assume the stability condition
$\rho:= {\mathbf E} \sigma_1/{\mathbf E} t_1 < m$ to hold. 

Again let $D$ be the stationary waiting time. We formulate the
following questions: what may the asymptotics for 
${\mathbf P} (D>x)$ be when $x$ is large and
what is the ``typical'' sample path which lead to such a large value
of the stationary waiting time. To answer (only partially!) these
questions, we need further restrictions on the distribution
of service times. 

First, we assume that the common distribution of service times is
heavy-tailed, i.e. for any $c>0$, its $c$th exponential moment does not 
exist, ${\mathbf E} e^{c\sigma_1}=\infty$.  
Whitt [\ref{WW}] formulated the following conjecture:
if $k \le \rho < k+1$  for an integer $k<m$, then
\begin{equation}\label{jump}
{\mathbf P} (D>x) \sim \gamma 
\left(
\overline{B}_I(\eta x) \right)^{m-k}
\quad \mbox{as} \quad x\to\infty
\end{equation}
``where $\gamma$ and $\eta$ are positive constants
(as functions of $x$)'' [sic, [\ref{WW}]] and where $\overline{B}_I$
is the tail of the integrated service time distribution.
Intuitively, formula (\ref{jump}) says that the main cause for
the stationary waiting time to be large is to have $m-k$ big
service times in the past.

In [\ref{FKor0}] and [\ref{FKor}] we show that if the distribution
of service times is {\it intermediate regularly varying} and if
$\rho$ is not an integer, then the conjecture of Whitt is correct with
$\gamma = \gamma (x)$ squeezed between two positive constants and
with $\eta$ being a constant. Also, the conjecture holds
if $\rho \in (0,1)$ and if the distribution $B_I$ is
any {\it subexponential} distribution. Also, we found that if $m=2$, $\rho<2$, $\rho \ne 1$,
and the service times distribution is again intermediately regularly
varying, then $\gamma$ is a constant.

The {\bf open problems} here are (I formulate them in the particular case of a
two server queue, $m=2$):
find the asymptotics for ${\mathbf P} (D>x)$  \\
(i) if $\rho=1$ -- at least, for 
some particular subexponential (say, regularly varying) distribution
of service times;\\
(ii) if $\rho \in (1,2)$ and if the distribution of service times
is heavy-tailed but has all power moments finite, for example,
if ${\mathbf P} (\sigma_1>x) = e^{-x^{\beta}}$, for some $\beta \in
(0,1)$.

It would be great to understand what are in these cases the ``typical'' 
paths which lead to large values of $D$.

\section{Further problems on multi-server queues}

Consider again the multi-server (say, 2-server) queue with
stationary and ergodic input $\{t_n, \sigma_n\}$, but
assume now that the discipline is
``join-the-shortest-queue'': there are individual queues in
front of the servers, and each arriving customer joins 
immediately the shortest queue (or one of the shortest at
random if there are many).

So, here are the {\bf open problems}: \\
-- how many stationary regimes
may exist if we do not assume
 any extra condition in addition to the obvious stability condition?\\
-- what are the minimal requirements for the uniqueness of
the stationary regime?\\
-- under what conditions (none?) do we have weak (or TV) convergence?

The model exhibits NO monotonicity.
It is not amenable to Loynes-type schemes.
It is entirely open ([\ref{TK}]).

\section{Greedy service mechanism}

There are many circumstances in our life where we may ask the
question, is a ``locally optimal'' (``greedy'') mechanisms
also optimal in the long-run? Below are examples of mathematical
models where the answer to such a question is open.


There are two continuous state space models where the stability conjecture is
obvious, but nobody is able to verify it. 
In both models, the driving algorithm contains a ``locally optimal''
(``greedy'') element. 
It looks like none of the existing 
stability
methods works here.

\subsection{Stability of a greedy server}

A single server is located on the circle. Particles arrive 
in a Poisson
stream of rate $\lambda$ and are uniformly distributed (as material points) 
on the
circle (people say that there is a ``Poisson rain'' of particles). 
It takes a single unit of time to serve a particle. After any service, the
particle disappears, the
server chooses to serve next the {\it closest} particle and moves 
to it with a (positive finite)
constant speed (ignoring new arrivals), serves it during
another unit of time, then chooses the next closest particle and moves to it, 
etc.

The {\bf conjecture} is: this model is stable for any $\lambda < 1$. A plausible
``proof'' might be as follows: if the number of requests is very large, then
the server is busy with service almost all the time (with a service speed close
to one), and then we may apply, say, fluid approximation ideas to deduce
the stability.

This model and this conjecture have already been known for more than 20 years,
see [\ref{CofGil}], 
but nobody has been able to succeed with obtaining either a proof or a counter-example
here. The key problem is the continuity of the state space, and there are
several results (see, e.g., [\ref{FL1}, \ref{FL2}, \ref{Scha}] 
for further details) with the proof
of a similar hypothesis for models with a finite state space (for instance,
you may replace the continuous circle by a finite lattice on it). If the
server uses any ``state-independent'' algorithm for moving (say, always
walks in the left direction or chooses the next direction with probability
$1/2$ independently of everything else), then it is easy to verify the
conjecture using the ideas explained above -- see, e.g., [\ref{FL}, \ref{KS}].

\subsection{Stability of a model with two streams:
stream of customers and stream of servers}

Again, there is a circle, but this time no server or service. 
Instead, there are two independent Poisson streams/rains, of ``black'' and 
of ``white'' particles, with rates $\lambda$ and $1$, respectively.
Black particles arrive at the circle and stop there, but white particles 
pass straight through the circle (this means they ``arrive and immediately 
disappear''). 
There is given a distance $\varepsilon >0$.
When a white particle passes through the circle at some point, it
observes all blacks in the $\varepsilon$-neighbourhood and takes
(deletes) the one which is the closest to itself (if there are any black particles at that instant).

The natural {\bf conjecture} is: stability should be guaranteed by the
condition $\lambda < 1$,  independently of the circle length and 
the number $\varepsilon$. But the problem is open too. Again, there exist
simple proofs for stability if the model is modified:
 if either the continuous state space (the circle) is replaces by
a finite set, or the greedy mechanism is replaced by any state-independent
mechanism (for instance, if a white particle takes one of blacks
from the neighbourhood ``at random'', with equal probabilities) 
 (see [\ref{VA}]).

\section{Stability may depend on the whole distribution}

The {\bf conjecture} is: even in simple queueing systems, 
the stability conditions
may depend both on the initial values and on the whole distribution
of the driving sequences.

Here is an example where the conjecture may be true (see [\ref{FCh}] for
more detail).

Consider a system with three servers (numbered 1 to 3) 
fed by a Poisson process with intensity
$\lambda$. There are three classes of customers and each arriving customer
becomes a class $i$ customer ($i=1,2,3$) with probability $1/3$.
Class $1$ customers may be served by 1st and 2nd servers (where 1st server
is ``left'' and 2nd server is ``right''), class 2 customers
by 2nd and 3rd servers (here 2nd server is ``left'' and 3rd is ``right''), 
and class 3 customers by 3d and 1st servers (here 3rd is ``left'' and 1st
is ``right''). 
Upon arrival, a customer chooses an accessible server with the shorter workload.
There are two probability distributions, $F_l$ and $F_r$, and a customer's
service time has distribution $F_l$ if it is served by 
its left server and $F_r$ otherwise. 

Simulations show that the conjecture may hold, but there is no a rigorous proof.

\section{Random fluid limits and positive Lebesgue measure
of the area of null-recurrence}

Consider an open polling system with two stations and 
two ``heterogeneous'' servers.
Each station $i=1,2$ has a Poisson input with intensity $\lambda_i=1$. 
For $i,j \in \{1,2\}$, service times of server $j$ at station $i$
are i.i.d. exponential with intensity 
$\mu_i^{(j)}$. Both servers follow the so-called exhaustive service 
policy:    
after completing a service, a server either starts with a service 
of a new customer
(if there is any), or leaves the station for the other one. We assume 
for simplicity
that ``walking'' (``switchover'') times are equal to zero. If there is no free
customer at either station, the server becomes ``passive''.

The system described has a nice fluid model where all fluid limits are random
and piecewise deterministic. The system is characterised by 4 parameters
$\{ \mu_i^{(j)}\}_{i,j=1,2}$, and
all of them have to be less than one in order to make the system stable (but
this is definitely not sufficient for stability). 

The {\bf conjecture} here is: in the positive 4-dimensional cube, the set of 
parameters $\{ \mu_i^{(j)}\}$ for which the system is ``null-recurrent''
has a positive Lebesgue measure. See [\ref{FKov}] for more details.

\section{Multi-access channel with protocols based on
partial information}

Consider a single channel which is shared among many users and transmits
packets (messages) of a single length. Assume that time is slotted
and each service time is equal to the slot length. The number of packages 
arriving into the system during a time slot is Poisson with parameter
$\lambda$. At the beginning of time slot $n$, each package is trying
to be transmitted with the same probability $p_n$ independently of
everything else. If two or more packages try to transmit simultaneously,
the transmissions collide and packages stay in the system and have
to try again later. If there is only one transmission, then
it is successful and the package leaves the system. If there are no
transmissions, then the slot is empty. Denote by
$W_n$ the number of packages in the system at the beginning of $n$th 
time slot. Assume that probabilities $p_n$ are defined inductively
and that $p_{n+1}$ may depend only on $p_n$ and the ``binary'' information of
whether there is a successful transmission in slot $n$ or not. 
Then the pairs $(W_n,p_n)$ form a time homogeneous Markov chain.

The {\bf open question} is: in the class of procedures (protocols)
described above, does there exist a protocol which makes the
Markov chain $(W_n,p_n)$ positive recurrent (ergodic), for
some positive $\lambda$? See [\ref{FT}] for more detail.

It is known [\ref{Hay}, \ref{Mikh}] that if a choice of $p_{n+1}$
is based on $p_n$ and of either of two other ``binary'' informations
($n$th slot was empty or not; or there was collision in $n$th slot or not),
then $\lambda < e^{-1}$ is necessary and sufficient for the existence
of an ergodic protocol.

\section*{\normalsize References}

\newcounter{bibcoun}
\begin{list}{\arabic{bibcoun}.}{\usecounter{bibcoun}\itemsep=0pt}
\small

\item\label{VA} V.~Anantharam, private communication.


\item\label{Bra} A.~Brandt.
On stationary waiting times and limiting behaviour of queues with many 
servers II. The $G/GI/m/\infty$ case. 
{\it Elektron. Informationsverarb. Kybernetik}, {\bf 21} (1985), 151--162.

\item\label{BFL} A.~Brandt, P.~Franken and B.~Lisek. 
{\it Stationary Stochastic Models}. 
Akademie-Verlag and Wiley, 1990.


\item\label{Bor} A.A.~Borovkov, {\it Asymptotic Methods in Queueing
Theory}, Nauka, 1976 (in Russian), Wiley, 1980  (English translation).

\item\label{CofGil} E.G.~Coffman and E.N.~Gilbert.
Polling and greedy servers on a line.
{\it Queueing Systems}, {\bf 2} (1987), 115--145.


\item\label{FL} G.~Fayolle and J.-M.~Lasgouttes. 
A state-dependent polling model with Markovian routing. 
In Frank P. Kelly and Ruth R. Williams, editors, 
IMA Volume 71 on Stochastic Networks, Springer (1995),
283--312.


\item\label{F83}  S.~Foss.
On Ergodicity Conditions for Multi-Server Queueing Systems.  
{\it Siberian Math. J.}, {\bf 24} (1983), 168--175. 

\item\label{FCh}
S.~Foss and N.~Chernova.
On stability of a partially accessible multi-station queue with 
state-dependent routing.  {\it Queueing Systems}, 
{\bf 29} (1998), 55--73. 

\item\label{FL1}
S.~Foss and G.~Last.
Stability of Polling Systems with State Dependent Routing and 
with Exhaustive Service Policies. 
{\it Annals of Applied Probability}, {\bf 6} (1996), 116--137.

\item\label{FL2} S.~Foss and G.~Last.
Stability of polling systems with general service policies and with 
state dependent routing. 
{\it Probab. Eng. Inform. Sci.}, {\bf 12} (1998), 49--68.

\item\label{FKon} S.~Foss and T.~Konstantopoulos.
An overview of some stochastic stability methods. 
{\it J. Oper. Res. Soc. Japan}, 
{\bf 47} (2004), 275--292.

\item\label{FKor0} S.~Foss and D.~Korshunov. 
Heavy tails in multi-server queues.
{\it Queueing Systems}, {\bf 52} (2006), 31--48. 

\item\label{FKor} S.~Foss and D.~Korshunov.
How big queues occur in a multi-server system with heavy tails.
(submitted)

\item\label{FKov} S.~Foss and A.~Kovalevskii. 
A stability criterion via fluid limits and its application to 
a polling model. {\it Queueing Systems}, {\bf 32} (1999), 131--168.

\item\label{FT} S.~Foss and A.~Tyurlikov, in preparation.

\item\label{KW}
 J.~Kiefer and J.~Wolfowitz. 
On the theory of queues with many servers.
{\it Trans. Amer. Math. Soc.}, {\bf 78} (1955), 1--18. 

\item\label{Hay} B.~Hajek.
Hitting-time and occupation-time bounds implied by
drift analysis with applications.
{\it Adv. Appl. Probab.}, {\bf 14} (1982), 502--525.

\item\label{TK} T.~Konstantopoulos, private communication.

\item\label{KS} D.~Kroese and V.~Schmidt.
 Queueing systems on a circle.
{\it J. Math. Methods Oper. Res.}
{\bf 37} (1993), 303--331.

\item\label{Loy} M.~Loynes. 
The stability of a queue with non- independent inter-arrival and service 
times. {\it Proc. Cambridge Phil. Soc.}, {\bf 58} (1962), 497--520. 

\item\label{MT} S.~Meyn and R.~Tweedie.
{\it Markov Chains and Stochastic stability.}
Springer, 1993.

\item\label{Mikh} W.~Mikhailov.
Geometric analysis of stability of Markov chains and its
applications.
{\it Probl. Inform. Transm.}, {\bf 24} (1988), 61--73.

\item\label{Naka} T.~Nakatsuka. 
The untraceable events method for absorbing processes.
{\it J. Appl. Prob.}, {\bf 43} (2006), 652--664.

\item\label{Scha} R.~Schassberger.
Stability of polling networks with state-dependent server routing.
{\it Probab. Eng. Inform. Sci.}, {\bf 9} (1995), 539--550.

\item\label{Sto} D.~Stoyan.
{\it Comparison methods for queues and other stochastic models}.
Springer, 1983.

\item\label{SW-V} A.~Scheller-Wolf and R.~Vesilo.
Structural interpretation and derivation of necessary and sufficient conditions
for delay moments in FIFO multiserver queues.
{\it Queueing Systems}, {\bf 54} (2006), 221--232.

\item\label{WW} W.~Whitt.
The impact of a heavy-tailed service-time distribution upon the
$M/GI/s$ waiting-time distribution.
{\it Queueing Systems}, {\bf 36} (2000), 71--87.

\end{list}


\end{document}